\documentclass[letterpaper,english]{article}
\usepackage[T1]{fontenc}
\usepackage[latin9]{inputenc}
\usepackage{amsmath}
\usepackage{amssymb}

\makeatletter

\pdfpageheight\paperheight
\pdfpagewidth\paperwidth

\makeatother

\usepackage{babel}
\begin{document}

\title{Faulhaber and Bernoulli\thanks{This work is licensed under the CC BY 4.0, a Creative Commons Attribution
License.}}

\author{Ryan Zielinski\\
ryan\_zielinski@fastmail.com}

\date{7 December 2018}
\maketitle
\begin{abstract}
In this note we will use Faulhaber's Formula to explain why the odd
Bernoulli numbers are equal to zero.
\end{abstract}

\section{Introduction}

For odd numbers greater than or equal to seven, why are the Bernoulli
numbers equal to zero? Because Faulhaber's Formula tells us that $\sum_{k=1}^{n}k^{2m+1}$
is a polynomial in $\left(\sum_{k=1}^{n}k\right)^{2}$, and $\left(\sum_{k=1}^{n}k\right)^{2}=\frac{n^{2}+2n^{3}+n^{4}}{4}$.

\section{Faulhaber's Formula}

We might know already 
\[
\sum_{k=1}^{n}k^{3}=\left(\frac{n(n+1)}{2}\right)^{2}=\left(\sum_{k=1}^{n}k\right)^{2}.
\]

\pagebreak{}

Through inductive reasoning like that in \cite{key-7} we might discover
further 
\begin{align}
2^{1}\cdot\left(\frac{n(n+1)}{2}\right)^{2} & =\binom{2}{1}\cdot\sum k^{3}\nonumber \\
2^{2}\cdot\left(\frac{n(n+1)}{2}\right)^{3} & =\binom{3}{0}\cdot\sum k^{3}+\binom{3}{2}\cdot\sum k^{5}\nonumber \\
2^{3}\cdot\left(\frac{n(n+1)}{2}\right)^{4} & =\binom{4}{1}\cdot\sum k^{5}+\binom{4}{3}\cdot\sum k^{7}\label{eq:1}\\
2^{4}\cdot\left(\frac{n(n+1)}{2}\right)^{5} & =\binom{5}{0}\cdot\sum k^{5}+\binom{5}{2}\cdot\sum k^{7}+\binom{5}{4}\cdot\sum k^{9}\nonumber \\
2^{5}\cdot\left(\frac{n(n+1)}{2}\right)^{6} & =\binom{6}{1}\cdot\sum k^{7}+\binom{6}{3}\cdot\sum k^{9}+\binom{6}{5}\cdot\sum k^{11}.\nonumber 
\end{align}
 (We abbreviate $\sum_{k=1}^{n}k^{m}$ by $\sum k^{m}$.) The general
case is 
\begin{eqnarray}
2^{m-1}\cdot\left(\frac{n(n+1)}{2}\right)^{m} & = & \binom{m}{0}\cdot\sum k^{m}+\binom{m}{2}\cdot\sum k^{m+2}\label{eq:2}\\
 &  & +\cdots+\binom{m}{m-3}\cdot\sum k^{2m-3}+\binom{m}{m-1}\cdot\sum k^{2m-1}\nonumber 
\end{eqnarray}
 or 
\begin{eqnarray}
2^{m-1}\cdot\left(\frac{n(n+1)}{2}\right)^{m} & = & \binom{m}{1}\cdot\sum k^{m+1}+\binom{m}{3}\cdot\sum k^{m+3}\label{eq:3}\\
 &  & +\cdots+\binom{m}{m-3}\cdot\sum k^{2m-3}+\binom{m}{m-1}\cdot\sum k^{2m-1}.\nonumber 
\end{eqnarray}

If we wish to prove such expressions, by \cite{key-1,key-2,key-3}
we may proceed using Pascal's observation of telescoping sums. Consider
the special case of 
\begin{eqnarray*}
\left(\frac{3\cdot4}{2}\right)^{m} & = & \left(\frac{1\cdot2}{2}\right)^{m}-0+\left(\frac{2\cdot3}{2}\right)^{m}-\left(\frac{1\cdot2}{2}\right)^{m}+\left(\frac{3\cdot4}{2}\right)^{m}-\left(\frac{2\cdot3}{2}\right)^{m}\\
 &  & =\sum_{k=1}^{3}\left[\left(\frac{k(k+1)}{2}\right)^{m}-\left(\frac{(k-1)k}{2}\right)^{m}\right].
\end{eqnarray*}
 For the general case of 
\begin{equation}
\left(\frac{n(n+1)}{2}\right)^{m}=\sum_{k=1}^{n}\left[\left(\frac{k(k+1)}{2}\right)^{m}-\left(\frac{(k-1)k}{2}\right)^{m}\right],\label{eq:4}
\end{equation}
 the expression in brackets we may rewrite as 
\[
\left(\frac{k}{2}\right)^{m}\left((1+k)^{m}-(-1+k)^{m}\right)
\]
 and then use the binomial theorem to arrive at (\ref{eq:2}) or (\ref{eq:3}).

\bigskip{}

What if we want to rewrite such expressions in terms of a particular
$\sum k^{2m+1}$? For example, by looking at (\ref{eq:1}) we may
write 
\begin{eqnarray*}
\sum k^{5} & = & \frac{1}{\binom{3}{2}}\cdot\left[2^{2}\cdot\left(\frac{n(n+1)}{2}\right)^{3}-\binom{3}{0}\cdot\sum k^{3}\right]\\
 &  & =\frac{1}{3}\cdot\left[4\cdot\frac{n(n+1)}{2}-1\right]\cdot\left(\sum k\right)^{2},
\end{eqnarray*}
 which then implies 
\begin{eqnarray*}
\sum k^{7} & = & \frac{1}{\binom{4}{3}}\cdot\left[2^{3}\cdot\left(\frac{n(n+1)}{2}\right)^{4}-\binom{4}{1}\cdot\sum k^{5}\right]\\
 &  & =\frac{1}{4}\cdot\left[2^{3}\cdot\left(\frac{n(n+1)}{2}\right)^{2}-4\cdot\frac{1}{3}\cdot\left(4\cdot\frac{n(n+1)}{2}-1\right)\right]\cdot\left(\sum k\right)^{2}.
\end{eqnarray*}
 By a proof by mathematical induction on the $m$ of (\ref{eq:2})
or (\ref{eq:3}) we may arrive at the general result of 
\begin{eqnarray}
\sum k^{2m+1} & = & \frac{1}{m+1}\cdot\left[2^{m}\cdot\left(\frac{n(n+1)}{2}\right)^{m-1}-a_{2}\cdot\left(\frac{n(n+1)}{2}\right)^{m-2}+a_{3}\cdot\left(\frac{n(n+1)}{2}\right)^{m-3}\right.\nonumber \\
 &  & \left.\mp\cdots\mp a_{m-2}\cdot\left(\frac{n(n+1)}{2}\right)^{2}\mp a_{m-1}\cdot\frac{1}{3}\cdot\left(4\cdot\frac{n(n+1)}{2}-1\right)\right]\cdot\left(\sum k\right)^{2},\label{eq:ff}
\end{eqnarray}
 where the $a_{i}$ are rational numbers and $m\geqq3$. This relationship
we will call Faulhaber's Formula. (For some of the history of the
problem, see \cite{key-1,key-2,key-3,key-4}.)

\section{Bernoulli Numbers}

By \cite{key-1,key-2} or Chapter 1 of \cite{key-6} we may define
the Bernoulli numbers $B_{n}$ by 
\[
B_{0}=1,\quad\sum_{k=0}^{n}\binom{n+1}{k}\cdot B_{k}=0,
\]
 where $n\geqq1$. For example, to find $B_{1}$ we write 
\[
\sum_{k=0}^{1}\binom{1+1}{k}\cdot B_{k}=\binom{2}{0}\cdot B_{0}+\binom{2}{1}\cdot B_{1}=0,
\]
 which implies $B_{1}=-\frac{1}{2}$. It turns out $B_{3}=B_{5}=0$.
Following the section ``Back to Faulhaber's form'' of \cite{key-4}
we may write 
\begin{eqnarray}
\sum k^{2m+1} & = & \frac{1}{2m+2}\cdot\left[\binom{2m+2}{0}\cdot B_{0}\cdot n^{2m+2}-\binom{2m+2}{1}\cdot B_{1}\cdot n^{2m+1}\right.\nonumber \\
 &  & \left.\pm\cdots+\binom{2m+2}{2m}\cdot B_{2m}\cdot n^{2}-\binom{2m+2}{2m+1}\cdot B_{2m+1}\cdot n\right],\label{eq:bn}
\end{eqnarray}
 for which we will assume $m\geqq3$.

\section{Conclusion}

With regard to the claim at the start, suppose we set (\ref{eq:ff})
and (\ref{eq:bn}) equal to one another. If we multiply out (\ref{eq:ff}),
we see it does not contain the term $n$. That means the last term
of (\ref{eq:bn}), 
\[
-\frac{1}{2m+2}\cdot\binom{2m+2}{2m+1}\cdot B_{2m+1}\cdot n,
\]
 must be equal to zero. In other words, for all $m\geqq3$, $B_{2m+1}=0$.

\section*{Acknowledgements}

The author thanks the anonymous referee for suggestions which improved
the quality of the paper.

\end{document}